\documentclass[12pt,leqno,fleqn]{amsart}  
\usepackage{amsmath,amstext,amsthm,amssymb,amsxtra}

\usepackage{txfonts} % pxfonts txfonts 
\usepackage[T1]{fontenc}
\usepackage{lmodern}

 \usepackage{euler}   % better than the option below

\usepackage{mathtools}
\mathtoolsset{showonlyrefs,showmanualtags}

\usepackage{hyperref} %,hypertexnames=false,colorlinks,[pagebackref]
\hypersetup{
    colorlinks=true,       % false: boxed links; true: colored links
    linkcolor=blue,          % color of internal links
    citecolor=magenta,        % color of links to bibliography
    filecolor=magenta,      % color of file links
    urlcolor=cyan           % color of external links
}

\usepackage[alphabetic,msc-links]{amsrefs}

\allowdisplaybreaks
\swapnumbers

\theoremstyle{plain} % definition 
\newtheorem{lemma}[equation]{Lemma} 
 
\newtheorem{theorem}[equation]{Theorem}

\theoremstyle{definition}
\newtheorem{definition}[equation]{Definition} 

\theoremstyle{remark}
\newtheorem{remark}[equation]{Remark}

\numberwithin{equation}{section}

\title {On the $ A_2$ inequality for Calder\'on-Zygmund Operators}
\author{Michael T. Lacey}   %  can use \and  

\address{ School of Mathematics, Georgia Institute of Technology, Atlanta GA 30332, USA}
\email {lacey@math.gatech.edu}
\thanks{Research supported in part by grant NSF-DMS 0968499.}

\begin{document}
\begin{abstract}
We prove  that for an $ L ^2 (\mathbb R ^{d}) $-bounded Calder\'on-Zygmund operator and  weight $ w\in A_2$, that we have the inequality 
below due to Hyt\"onen,   
\begin{equation*}
\lVert T  \rVert_{ L ^{2} (w) \to L ^{2} (w)} 
\le C _{T} [w]_{A_2}   \,. 
\end{equation*}
Our proof will appeal to a distributional inequality used by several authors,  adapted Haar functions, and standard 
stopping times.  
\end{abstract}
\maketitle 

%%%%%%%%%%%%%%%%%%%%%%%%%%%%%% SECTION  SECTION SECTION
%%%%%%%%%%%%%%%%%%%%%%%%%%%%%% SECTION  SECTION SECTION 
\section{Introduction: Main Theorem} %\label{s.}

We are interested in estimates for the norms of Calder\'on-Zygmund operators on weighted $ L ^{p}$-spaces, 
a question that has attracted significant interest recently;  definitive estimates 
of this type were first obtained in   \cites{1007.4330}, with a range of prior and subsequent contributions.  
In this paper, we will concentrate on $ p=2$, and give a new proof, more elementary than some of the preceding proofs. 

 Let  $ w$ be a weight on $ \mathbb R ^{d}$  with density also written as $ w$.  Assume $ w>0$ a.\thinspace e. We define  $ \sigma= w ^{-1}  $, which is defined a.\thinspace e.\thinspace, and set 
\begin{equation} \label{e.A2}
[w]_{A_2} := \sup _{Q} \frac {w (Q)} {\lvert  Q\rvert }  \frac {\sigma (Q)} {\lvert  Q\rvert } \,. 
\end{equation}

We give a new proof of 

%%%%%%%%%%%%%%%%%%%%%%%%%%%%%% THEOREM THEOREM THEOREM
\begin{theorem}\label{t.main}[\cite{1007.4330}]
Let  $  T $ be an $ L ^2 $ bounded Calder\'on-Zygmund operator, and  $ w\in A_2$.   It then holds that 
\begin{align}
\lVert  T  f  \rVert _{ L ^{2} (w)} &\le C_{T} [w ]_{ A _{2}}
\lVert  f \rVert _{L ^{2} (w)}.
\end{align} 
\end{theorem}
%%%%%%%%%%%%%%%%%%%%%%%%%%%%%% THEOREM THEOREM THEOREM

All proofs in this level of generality have used Hyt\"onen's 
random Haar shift representation from \cite{1007.4330}. So does this proof.  After this point, two strategies of prior proofs 
are  (a)  fundamental appeal to two-weight inequalities, an approach 
initiated in \cite{MR2657437} and further refined in \cite{1007.4330,1103.5229,1103.5562,1010.0755}, 
or (b) constructions of appropriate Bellman functions \cite{1105.2252}, extending the works of \cites{MR2354322,MR2433959}.  
In our approach, we borrow the distributional inequalities central to the two-weight approach, but then combine them 
with adapted Haar functions from the Bellman approach. Then, the familiar stopping time  considerations of \S\ref{s.proof} 
are sufficient to  conclude the proof.  More detailed  histories of this question can be found in the introductions to \cite{1010.0755,1105.2252,1103.5562}.  

As we will concentrate on the case of  $ L ^2 $ estimates, we will frequently use the notation $ \lVert f\rVert_{w} 
:= \Bigl[ \int f ^2 \; w (dx) \Bigr] ^{1/2} $.  At one or two points,  an $ L ^{1}$ norm is needed, and this will be 
clearly indicated.

%%%%%%%%%%%%%%%%%%%%%%%%%%%%%% SECTION  SECTION SECTION
%%%%%%%%%%%%%%%%%%%%%%%%%%%%%% SECTION  SECTION SECTION 
\section{Haar Shift Operators} \label{s.haarShift}

In this section, we introduce fundamental dyadic approximations of Calder\'on--Zygmund operators, the Haar shifts, and state reduction of the Main Theorem~\ref{t.main} to a similar statement, Theorem~\ref{t.haarShiftWtd},   in this dyadic model. In so doing, we are following 
the lead of \cite{1007.4330}.

%%%%%%%%%%%%%%%%%%%%%%%%%%%%%%  DEFINITION DEFINITION DEFINITION
\begin{definition}\label{d.grid} A \emph{dyadic grid} is a collection $ \mathcal D$ of cubes so that for each $Q$ we have 
that 
%%  ENUMERATE
\begin{enumerate}
\item  The set of cubes $ \{Q'  \in \mathcal D\;:\;  \lvert  Q'\rvert= \lvert  Q\rvert  \}$ partition $ \mathbb R ^{d}$, ignoring overlapping 
boundaries of cubes. 
\item $ Q$ is a  union of  cubes in a collection $ \textup{Child} (Q)  \subset \mathcal D$, called the \emph{children of $ Q$}. 
There are $ 2 ^{d}$ children of $ Q$, each of  volume $ \lvert  Q'\rvert= 2 ^{-d} \lvert  Q\rvert  $. 
\end{enumerate}
%% ENUMERATE
We refer to any subset of a dyadic grid as simply a \emph{grid}.  
\end{definition}
%%%%%%%%%%%%%%%%%%%%%%%%%%%%%%  DEFINITION DEFINITION DEFINITION

The standard choice for $ \mathcal D$ consists of the cubes $ 2 ^{k} \prod _{s=1} ^{d} [n_s, n_s+1)$ for $ k, n_1 ,\dotsc, n_d\in \mathbb Z $.  
But, the reduction we are stating here  depends upon a random family of dyadic grids.    This next definition is at slight variance with 
that of \cites{1103.5229,1010.0755,1007.4330}.  

%%%%%%%%%%%%%%%%%%%%%%%%%%%%%%  DEFINITION DEFINITION DEFINITION
\begin{definition}\label{d.haarShift} For integers $(m,n) \in \mathbb Z _+ ^2 $, we say that  a linear operator $ \mathbb S $ is a \emph{(generalized) Haar shift operator of   complexity type $ (m,n)$} if 
\begin{equation}\label{e.mn}
\begin{split}
  \mathbb S  f (x) &= \sum_{Q \in \mathcal D}\mathbb{S}_Q f(x)
= 
 \sum_{Q \in \mathcal D}\int_{Q}s_{Q}\left( x,y\right) f\left( y\right) dy 
\end{split}
\end{equation}
where  here and throughout $\ell (Q)= \lvert  Q\rvert ^{1/d} $, and   these properties hold. 
\begin{enumerate}
  \item  $s_{Q}$, the kernel of the component $\mathbb{S}_Q$, is supported on $Q\times Q$ and $\Vert s_{Q}\Vert
_{\infty }\leq \frac{1}{\left\vert
Q\right\vert }$. It is easy to check that
\begin{equation*} 
\Biggl\lvert  \sum_{Q \in \mathcal D}\int_{Q}s_{Q}\left( x,y\right) \Biggr\rvert\lesssim\frac{1}{|x-y|^d}.
\end{equation*}
 \item The kernel  $ s_Q$ is constant on dyadic rectangles $ R \times S \subset Q \times Q $ with $ \ell (R) \le 2 ^{-m} \ell (Q)$ 
 and $ \ell (S) \le 2 ^{-n} \ell (Q)$.  
 \item For any subset $ \mathcal D'\subset \mathcal D$, it holds that we have 
 \begin{equation*}
\Bigl\lVert  \sum_{Q\in \mathcal D'} \mathbb S _{Q} f \Bigr\rVert_{2} \le \lVert f\rVert_{2} \,. 
\end{equation*}
\end{enumerate}
We say that the \emph{complexity} of $ \mathbb S $ is $  \kappa := 1+\max (m,n)$.  
\end{definition}
%%%%%%%%%%%%%%%%%%%%%%%%%%%%%%  DEFINITION DEFINITION DEFINITION

Note that the last property above is an statement about unconditionality of the sum in the operator norm. This is in fact 
a standard part of Calder\'on-Zygmund theory---and one that is automatic, depending upon exactly how the 
definition is formulated.  This property is fundamental to the proofs of this paper, and other results that we merely 
cite, justifying our inclusion of this property into the definition.

The main results of \cite{1010.0755} (see \cite{1010.0755}*{Theorem 4.1}; also \cite{1007.4330}*{Theorem 4.2}) allows us 
to reduce  the proof of the Main Theorem~\ref{t.main}  to the verification of the following dyadic variant.

%%%%%%%%%%%%%%%%%%%%%%%%%%%%%% THEOREM THEOREM THEOREM
\begin{theorem}\label{t.haarShiftWtd}  Let $ \mathbb S $ be a Haar shift operator with complexity $ \kappa $. 
For $ w \in A_2$,  we then have the estimates 
\begin{align} \label{e.haarShiftWtd}
\lVert \mathbb S f \rVert_{w} & \lesssim \kappa  [w] _{A_2} 
\lVert f\rVert_{w} 
\end{align}
\end{theorem}
%%%%%%%%%%%%%%%%%%%%%%%%%%%%%% THEOREM THEOREM THEOREM

Indeed, any polynomial dependence on the complexity parameter $\kappa$ would suffice for Theorem~\ref{t.main}.
(The linear bound in $ \kappa $ was shown in \cite{1103.5229}*{Theorem 2.10} in even greater generality in $ L ^{p}$ and maximal 
truncations. Later, and by different methods, it was 
 shown by \cite{1105.2252} as stated above.) 

In the remainder of this paper, $ \mathbb S $ will denote a Haar shift operator of complexity $ \kappa $, 
with \emph{scales separated by $ \kappa $}.  Namely, we have for a subset $ \mathcal D _{\kappa } \subset \mathcal D$, 
\begin{equation} \label{e.separated}
\mathbb{S}f\left( x\right) =\sum_{Q\in \mathcal{D} _{\kappa }}\int_{Q}s_{Q}\left( x,y\right) f\left( y\right) dy,
\end{equation}
and $  \mathcal D _{\kappa }$ consists of all dyadic intervals with $ log _{2} \ell (Q) = \ell \mod \kappa $, for some fixed 
integer $ 0\le \ell < \kappa $. In particular, if $ Q', Q \in \mathcal D _{\kappa }$ and $ Q'\subsetneqq Q$, then $ \int _{Q} s _{Q} (x,y) f (y)$ 
is constant on $ Q' $. The dual statement is also true.

%%%%%%%%%%%%%%%%%%%%%%%%%%%%%% SECTION  SECTION SECTION
%%%%%%%%%%%%%%%%%%%%%%%%%%%%%% SECTION  SECTION SECTION 
\section{The Basic Inequalities} \label{s.dist}
We make a remark here about the formulation of the inequalities that we will consider below.  Recalling the dual weight 
$ \sigma = w ^{-1} $ to an $ A_2$ weight, we will show that 
\begin{equation} \label{e.dual}
\lVert \mathbb S (f \sigma ) \rVert_{w} \lesssim \kappa  [w] _{A_2} 
\lVert f\rVert_{\sigma }  
\end{equation}
This is formally equivalent to the statement we are proving, namely \eqref{e.haarShiftWtd}, moreover the inequality above is 
the natural way to phrase the inequality as it dualizes in the natural way: Interchange the roles of $ w$ and $ \sigma $.  
Accordingly, we will especially in the next section, use the notation $ \langle f , g \rangle _{w}$ for the natural inner-product on $ L ^2 (w)$.

The arguments initiated in \cite{MR2657437}, further refined in \cites{1007.4330,1010.0755,1103.5229}, yield the following 
estimates for Haar shifts on intervals.  

%%%%%%%%%%%%%%%%%%%%%%%%%%%%%% LEMMA LEMMA LEMMA
\begin{lemma}\label{l.I} Let  $ w\in A_2$, $ \mathbb S $ a Haar shift operator of complexity $ \kappa $  as in \eqref{e.separated}. 
 For a cube $ Q$, and let $ \mathcal Q \subset \mathcal D _{\kappa }$ 
be a collection of cubes contained in $ Q$. We have 
\begin{align}\label{e.LL1} 
\int _{Q} \lvert  \mathbb S _{\mathcal Q} (\mathbf 1_{Q}) \rvert\; w (dx) & \lesssim [w] _{A_2} \lvert  Q\rvert   \,, 
\\ \label{e.LL2}
\int _{Q}  \mathbb S _{\mathcal Q} (\mathbf 1_{Q}) ^2 \; w (dx) & \lesssim [w] _{A_2} ^2  \sigma (Q) \,. 
\end{align}
\end{lemma}
%%%%%%%%%%%%%%%%%%%%%%%%%%%%%% LEMMA LEMMA LEMMA

There are two estimates, one of an $ L ^{1} (w)$ norm, with the right hand side being Lebesgue measure.  The second is 
an $ L ^2 (w)$ norm, with the right hand side being in terms of $ \sigma $.  
Indeed, the papers \cite{MR2657437,1007.4330,1103.5229,1010.0755} are argued such that  the second estimate 
\eqref{e.LL2}, combined with a general two-weight theorem, imply the linear bound in $ A_2$.  Thus, the point of this paper is that 
the general two-weight theorems are not needed.   

We suppress the proof of the $ L ^2 $ estimate, which is based upon a corona decomposition, and distributional estimate. 
The $ L ^{1} $ estimate follows from the same line of attack.  
The reader can consult for instance \cite{1010.0755}*{Lemma 5.7} or \cite{1103.5229}*{Section 11}.  
Of course the $ L ^{2} $ estimate implies an estimate for the $ L ^{1} $ norm, and it is interesting to note that it is worse than what one gets by using the proof  of the $ L ^2 $ estimate.

%%%%%%%%%%%%%%%%%%%%%%%%%%%%%% SECTION  SECTION SECTION
%%%%%%%%%%%%%%%%%%%%%%%%%%%%%% SECTION  SECTION SECTION 
\section{Proof of the Weighted Estimate for the Haar Shift Operators} \label{s.proof}

We will need the martingale difference operators associated with $ \mathcal D _{\kappa }$, 
and weight $ \sigma $.   For cube $ Q\in \mathcal D_ \kappa $ 
consider the martingale difference operator 
\begin{equation*}
D ^{\sigma } _{Q} f 
:= \sum_{\substack{Q'\subset Q\\ \ell (Q') = 2 ^{- \kappa}\ell (Q) }}  
\mathbb E ^{\sigma } _{Q'} f  \mathbf 1_{Q'} - \mathbb E ^{\sigma } _{Q} f \,. 
\end{equation*}
Here, $ \mathbb E ^{\sigma } _{Q} f = \sigma (Q) ^{-1} \int _{Q} f \; \sigma dx $. 
The operators $ D ^{\sigma } _{Q} f$ are self-adjoint contractions on $ L ^{2} (\sigma )$, and  satisfy 
the fundamental orthogonality relationship is 
\begin{equation}\label{e.sumD}
\sum_{Q\in \mathcal D_ \kappa } \lVert  D ^{\sigma } _{Q} f \rVert_{\sigma } ^2 \le \lVert f\rVert_{\sigma } ^2 \,,  
\end{equation}
which holds under minimal assumptions on $ \sigma $, satisfied for a weight with density non-negative almost everywhere. 

Now, complexity shows that for any fixed $ Q$, the components of the Haar shift operator are 
\begin{equation*}
\int _{Q} s _{Q} (x,y) f (y) \; \sigma (dy) 
= \mathbb E _{Q} ^{\sigma } f 
\int _{Q} s _{Q} (x,y)  \; \sigma (dy)  +  
\int _{Q} s _{Q} (x,y) D ^{\sigma } _{Q} f (y) \; \sigma (dy) 
\end{equation*}
Note that 
the bilinear form $ \langle \mathbb S _{\sigma } f, g \rangle _{w}$  is the linear combination of  the 
three  terms below, and their duals.  
\begin{gather}\label{e.main}
\langle \mathbb U _{\sigma } f, g \rangle _{w} := 
\sum_{Q\in \mathcal D} \mathbb E _{Q} ^{\sigma } f 
\int _{Q}\int _{Q} s _{Q} (x,y)  g (x) \; \sigma (dy)  w (dy)
\\ \label{e.D}
\mathbb V _{\sigma} (f,g) 
:= 
\sum_{Q\in \mathcal D} \mathbb E _{Q} ^{\sigma } f 
\int _{Q}\int _{Q} s _{Q} (x,y)  D ^{w} _{Q}  g (x) \; \sigma (dy)  w (dy) 
\\ \label{e.B}
\mathbb W   (f,g) := 
\int _{Q} s _{Q} (x,y) D ^{\sigma } _{Q} f (y) D^{w} _{Q} g (x) \; \sigma (dy)  w (dx) 
\end{gather}
By dual, we mean that the roles of $ w$ and $ \sigma $ are interchanged, which is relevant to $ \mathbb U _{\sigma }$ 
and $ \mathbb V	 _{\sigma }$ above.  We will show that each of these three bilinear forms is bounded by $ [w] _{A_2} \lVert f\rVert_{\sigma } \lVert g\rVert_{w}$, which estimate also applies to their duals.  
Recall that we have assumed the separation of scales condition \eqref{e.separated}, hence under this condition we have proved 
Theorem~\ref{t.haarShiftWtd} with absolute constant. This proves the Theorem as stated, with $ \kappa $ bound 
in terms of complexity.

We prove the difficult estimate first, the estimate for $ \mathbb U  _{\sigma }$.  

%%%%%%%%%%%%%%%%%%%%%%%%%%%%%% SUBSECTION SUBSECTION SUBSECTION SUBSECTION
 %%%%%%%%%%%%%%%%%%%%%%%%%%%%%% SUBSECTION SUBSECTION SUBSECTION SUBSECTION 
\subsection*{The bound for $ \mathbb U _{\sigma }$.}%\label{ss.}

The essential tool is this  corona decomposition. 

%%%%%%%%%%%%%%%%%%%%%%%%%%%%%%  DEFINITION DEFINITION DEFINITION
\begin{definition}\label{d.fStopping} 
We say that $\mathcal  F \subset \mathcal D _{\kappa }$ is a set of $ f$-stopping cubes if 
these conditions are met.  
%%  ENUMERATE
\begin{enumerate}
\item  If $ F, F'\in \mathcal F$, $ F'\subsetneqq F$ then 
$ \rho (F) :=  \mathbb E _{F} ^{\sigma } \lvert  f\rvert  > 4 \mathbb E ^{\sigma }  _{Q_0} \lvert  f\rvert  $. 

\item Every cube $ Q\in \mathcal D _{\kappa }$ is contained in some $ F\in \mathcal F$. 

\item Let $ \mathcal D _F$ be those cubes for which $ F$ is the minimal element of $ \mathcal F$ containing $ Q$. For every 
$ Q\in \mathcal D_F$, we have $ \mathbb E ^{\sigma } _{Q} \lvert  f\rvert  \le 4 \mathbb E ^{\sigma } _{F} \lvert  f\rvert  $. 
\end{enumerate}
%% ENUMERATE
It is easy to recursively construct such a collection $\mathcal F $, for $ \sigma \in A_2$, which is the case we are considering.  
\end{definition}
%%%%%%%%%%%%%%%%%%%%%%%%%%%%%%  DEFINITION DEFINITION DEFINITION

A basic fact, a consequence of the  maximal function estimate for general weights, that we have  
\begin{equation}\label{e.F<}
\Bigl\lVert \sum_{F\in \mathcal F} \rho (F)  \mathbf 1_{F} \Bigr\rVert_{\sigma } ^2  \lesssim 
\sum_{F\in {\mathcal F}}  \rho  (F)^2 \sigma (F) \lesssim \lVert f\rVert_{\sigma } ^2\,, 
\qquad  \rho  (F) = \mathbb E ^{\sigma } _{F} \lvert  f\rvert ^2 \,. 
\end{equation}

The collections $ \mathcal D_F$ give a decomposition of $ \mathbb U _{\sigma }$ via   
\begin{align*}
\mathbb U _{\sigma , F} f &:= \sum_{Q\in \mathcal D_F} 
 \mathbb E _{Q} ^{\sigma } f 
\int _{Q} s _{Q} (x,y)\; \sigma (dy) 
\\
& = \rho (F) 
 \sum_{Q\in \mathcal D_F} 
[ \mathbb E _{Q} ^{\sigma } f  \cdot \rho (F) ^{-1} ] 
\int _{Q} s _{Q} (x,y)\; \sigma (dy) 
\,. 
\end{align*}
Note that the products $ \lvert  \mathbb E _{Q} ^{\sigma } f  \cdot \rho (F) ^{-1}\rvert   $ are never more than $ 4$, 
so by unconditionality  of Haar shifts, the integral estimates of the previous section apply to the expressions above.

We abandon duality,  expanding  
\begin{align*}
\lVert U _{\sigma } f\rVert_{w} ^2 & \le 
\Bigl\lVert  \sum_{F\in \mathcal F} \bigl\lvert\mathbb U _{\sigma , F} f  \bigr\rvert \Bigr\rVert_{w} ^2 
 \le I + 2 I\!I 
\\
I & :=  \sum_{F\in \mathcal F} \bigl\lVert\mathbb U _{\sigma , F} f  \bigr\rVert _{w} ^2  \,,
\\
 I\!I &:= 
 \sum_{F\in \mathcal F} \sum_{\substack{F'\in \mathcal F\\ F'\subsetneqq F }}\int _{F'} 
\bigl\lvert \mathbb U _{\sigma , F} f  \mathbb U _{\sigma , F'} f   \bigr\rvert  \; w (dx) \,. 
\end{align*}
These are the diagonal and off-diagonal terms.  The diagonal is immediate from \eqref{e.LL2} and \eqref{e.F<}: 
\begin{align*}
I &\lesssim 
[w] _{A_2} ^2 
 \sum_{F\in \mathcal F}  \rho (F) ^2  \sigma (F) \lesssim  [w] _{A_2} ^2  \lVert f\rVert_{w} ^2 \,. 
\end{align*}

The off-diagonal is as follows.  By the separation of scales hypothesis, 
note that in the definition of $  I\!I$, that $ \mathbb U _{\sigma , F} f  $ is constant on $ F'$ in the display below. 
Hence, by  \eqref{e.LL1}, we have 
\begin{align*}
 I\!I & \lesssim 
  [w] _{A_2}
 \sum_{F\in \mathcal F} \sum_{\substack{F'\in \mathcal F\\ F'\subsetneqq F }} 
 \mathbb E _{F'} ^{\sigma } 
 \bigl\lvert \mathbb U _{\sigma , F} f    \bigr\rvert \cdot 
 \rho (F')  \lvert  F'\rvert  
 \\
 & \lesssim   [w] _{A_2}
 \int \sum_{F\in \mathcal F}  \bigl\lvert \mathbb U _{\sigma , F} f    \bigr\rvert \cdot   \phi \; dx  
 \\
 \noalign{\noindent where $ \phi := \sum_{F\in \mathcal F} \rho(F)  \mathbf 1_{F}$, and using the identity $ w \cdot \sigma \equiv 1$,}
  & =   [w] _{A_2}
 \int \sum_{F\in \mathcal F}  \bigl\lvert \mathbb U _{\sigma , F} f    \bigr\rvert \cdot   \phi \; 
\sqrt {w (x) \sigma (x)}
 dx  
 \lesssim 
\Biggl\lVert  \sum_{F\in \mathcal F} \bigl\lvert\mathbb U _{\sigma , F} f  \bigr\rvert \Biggr\rVert_{w}
\lVert \phi \rVert_{\sigma } \,. 
\end{align*}
We have however $ \lVert \phi \rVert_{\sigma } \lesssim \lVert M ^{\sigma } f \rVert_{\sigma } \lesssim \lVert f\rVert_{\sigma }$.  
Combining estimates, we see that we have proved 
\begin{align*}
\Bigl\lVert  \sum_{F\in \mathcal F} \bigl\lvert\mathbb U _{\sigma , F} f  \bigr\rvert \Bigr\rVert_{w} ^2 
\lesssim 
  [w] _{A_2}  ^2 
\lVert f\rVert_{\sigma } ^2  
+ [w] _{A_2}  \Bigl\lVert  \sum_{F\in \mathcal F} \bigl\lvert\mathbb U _{\sigma , F} f  \bigr\rvert \Bigr\rVert_{w} \lVert f\rVert_{\sigma } 
\end{align*}
which immediately implies our linear bound in $ A_2$ for the 
term $ \mathbb U _{\sigma }$.

%%%%%%%%%%%%%%%%%%%%%%%%%%%%%% SUBSECTION SUBSECTION SUBSECTION SUBSECTION
 %%%%%%%%%%%%%%%%%%%%%%%%%%%%%% SUBSECTION SUBSECTION SUBSECTION SUBSECTION 
\subsection*{The Remaining Estimates}%\label{ss.}

%%%%%%%%%%%%%%%%%%%%%%%%%%%%%% PROOF PROOF PROOF
\begin{proof}[The bound for $ \mathbb V _{\sigma ,k} (f, g)$.]  We consider $ \mathbb V _{\sigma ,k} (f, g)$, defined in \eqref{e.D}.  
Using the orthogonality property of martingale differences \eqref{e.sumD}, we see that 
\begin{align*}
\bigl\lvert \mathbb V _{\sigma ,k} (f, g)\bigr\rvert 
& \le 
\sum_{Q\in \mathcal D} \lvert   \mathbb E _{Q} ^{\sigma } f \rvert \cdot 
\Bigl\lvert \bigl\langle       \int _{Q} s _{Q} (x,y) \; \sigma (dy)  
, D ^{w} _{Q}  g  \bigr\rangle _{w}\Bigr\rvert
\\
& \le 
\lVert g\rVert_{w} 
\Biggl\lVert 
\Biggl[
\sum_{Q\in \mathcal D}   \Bigl[ \mathbb E _{Q} ^{\sigma } f \cdot 
\int _{Q} s _{Q} (x,y) \; \sigma (dy)   \Bigr] ^2 
\Biggr] ^{1/2} 
\Biggr\rVert_{w} 
 \lesssim 
[w] _{A_2} \lVert f\rVert_{\sigma } \lVert g\rVert_{w} \,. 
\end{align*}
The last line follows from the bound already proved for the operator $ \mathbb U _{\sigma }$ and the 
 unconditionality. By a standard averaging over random choices of signs, 
we can deduce the linear in $ A_2$ bound for the square function above.  
\end{proof}
%%%%%%%%%%%%%%%%%%%%%%%%%%%%%% PROOF PROOF PROOF

%%%%%%%%%%%%%%%%%%%%%%%%%%%%%% PROOF PROOF PROOF
\begin{proof}[The bound for $ \mathbb W$.] We insert $ \sqrt {w \sigma }$ into the integrals below, and 
use the bound $ \lvert s _{Q} (x,y)\rvert \le \lvert  Q\rvert ^{-1}   $, to see that 
\begin{align*}
\Biggl\lvert 
\int _{Q} \int _{Q} s _{Q} (x,y) D ^{\sigma } _{Q} f (y) D^{ w} _{Q} g (x) \; dx dy
\Biggr\rvert & \le
\frac 1 {\lvert  Q\rvert } \int_Q \lvert  D ^{\sigma } _{Q} f (y) \rvert \; dy 
\int _{Q} \lvert  D^{ w} _{Q} g (x)\rvert  \; dx 
\\
& \le 
\lVert D ^{\sigma } _{Q} f  \rVert_{\sigma } 
\lVert D^{ w} _{Q} g \rVert_{w}
\Biggl[
\frac {\sigma (Q)} {\lvert  Q\rvert } \frac {w (Q)} {\lvert  Q\rvert }
\Biggr] ^{1/2} 
\\
& \le [w] _{A_2} \lVert D ^{\sigma } _{Q} f  \rVert_{\sigma } 
\lVert D^{ w} _{Q} g \rVert_{w} \,, 
\end{align*}
since we always have $ [w] _{A_2} \ge 1$.  
The martingale differences are pairwise orthogonal in $ L ^{2} (\sigma )$, and $ L ^{2} (w)$, so that a second 
application, in the variable $Q$, of the Cauchy-Schwartz inequality finishes this case. 
\end{proof}
%%%%%%%%%%%%%%%%%%%%%%%%%%%%%% PROOF PROOF PROOF

%%%%%%%%%%%%%%%%%%%%%%%%%%%%%% REMARK REMARK REMARK
\begin{remark}\label{r.oneOther} Rather than consider the operators $ \mathbb V _{\sigma }$ in \eqref{e.D}, and the dual expression, 
we could have considered 
\begin{equation*}
\widetilde {\mathbb V } (f,g) 
:= \sum_{Q\in \mathcal D _{k}} \mathbb E _{Q} ^{\sigma }   f \cdot 
 \int _{Q}\int _{Q} s _{Q} (x,y) \sigma (dy) w (dx) \cdot  \mathbb E _{Q} ^{w} g \,. 
\end{equation*}
It follows from unconditionality of Haar shift operators, and the estimate \eqref{e.LL1}, that we have the uniform estimate 
\begin{equation*}
\sum_{\substack{Q\in \mathcal D_ \kappa \\ Q \subset Q_0 }} 
\Biggl\lvert  \int _{Q}\int _{Q} s _{Q} (x,y) \sigma (dy) w (dx)   \Biggr\rvert 
\lesssim [w] _{A_2} \lvert  Q_0\rvert \,, \qquad \mathcal Q_0\in \mathcal D.  
\end{equation*}
From this, it is easy to see that 
\begin{align*}
\bigl\lvert \widetilde {\mathbb V } (f,g)  \bigr\rvert & \lesssim [w] _{A_2}
\int M ^{\sigma } f \cdot M ^{w} g \; dx 
\\
& =  [w] _{A_2}
\int M ^{\sigma } f \cdot M ^{w} g \; (w (x) \sigma (s)) ^{1/2}  dx 
\\
& \le  [w] _{A_2} \lVert M ^{\sigma }f\rVert_{\sigma } \lVert M ^{w}g\rVert_{w} 
 \lesssim [w] _{A_2} \lVert f\rVert_{\sigma } \lVert g\rVert_{w} \,. 
\end{align*}
Compare to Section 4 of \cite{1103.5347}.  But, we do not prefer this proof as it obscures the central role of the operator 
$ \mathbb U _{\sigma }$.  
\end{remark}

%%%%%%%%%%%%%%%%%%%%%%%%   references 


\begin{bibsection}
\begin{biblist}


\bib{MR2433959}{article}{
   author={Beznosova, Oleksandra V.},
   title={Linear bound for the dyadic paraproduct on weighted Lebesgue space
   $L\sb 2(w)$},
   journal={J. Funct. Anal.},
   volume={255},
   date={2008},
   number={4},
   pages={994--1007},
   issn={0022-1236},
   review={\MR{2433959}},
}

%\bib{MR1124164}{article}{
%  author={Buckley, Stephen M.},
%  title={Estimates for operator norms on weighted spaces and reverse Jensen inequalities},
%  journal={Trans. Amer. Math. Soc.},
%  volume={340},
%  date={1993},
%  number={1},
%  pages={253--272},
%  issn={0002-9947},
%}




\bib{MR2628851}{article}{
  author={Cruz-Uribe, David},
  author={Martell, Jos{\'e} Mar{\'{\i }}a},
  author={P{\'e}rez, Carlos},
  title={Sharp weighted estimates for approximating dyadic operators},
  journal={Electron. Res. Announc. Math. Sci.},
  volume={17},
  date={2010},
  pages={12--19},
  issn={1935-9179},
}


\bib{1001.4254}{article}{
   author={Cruz-Uribe, David},
   author={Martell, Jos{\'e} Mar{\'{\i}}a},
   author={P{\'e}rez, Carlos},
   title={Sharp weighted estimates for classical operators},
   date={2010},
      eprint={http://arxiv.org/abs/1001.4254},
}


%\bib{MR2281449}{article}{
%   author={Dragi{\v{c}}evi{\'c}, Oliver},
%   author={Petermichl, Stefanie},
%   author={Volberg, Alexander},
%   title={A rotation method which gives linear $L^p$ estimates for powers
%   of the Ahlfors-Beurling operator},
%   language={English, with English and French summaries},
%   journal={J. Math. Pures Appl. (9)},
%   volume={86},
%   date={2006},
%   number={6},
%   pages={492--509},
%   issn={0021-7824},
%   review={\MR{2281449 (2007k:30074)}},
%  %doi={10.1016/j.matpur.2006.10.005},
%}


%\bib{MR1992955}{article}{
%   author={Dragi{\v{c}}evi{\'c}, Oliver},
%   author={Volberg, Alexander},
%   title={Sharp estimate of the Ahlfors-Beurling operator via averaging
%   martingale transforms},
%   journal={Michigan Math. J.},
%   volume={51},
%   date={2003},
%   number={2},
%   pages={415--435},
%   issn={0026-2285},
%   review={\MR{1992955 (2004c:42030)}},
%  %doi={10.1307/mmj/1060013205},
%}


%
%\bib{MR727244}{article}{
%  author={Hru{\v {s}}{\v {c}}ev, Sergei V.},
%  title={A description of weights satisfying the {$A_{\infty }$} condition of {M}uckenhoupt},
%  journal={Proc. Amer. Math. Soc.},
%  fjournal={Proceedings of the American Mathematical Society},
%  volume={90},
%  year={1984},
%  number={2},
%  pages={253--257},
%  issn={0002-9939},
%  coden={PAMYAR},
%  mrclass={42B30},
%  review={\MR{727244 (85k:42049)}},
%  mrreviewer={R. Anantharaman},
%  %doi={10.2307/2045350},
%  url={http://dx.doi.org/10.2307/2045350},
%}

\bib{MR0312139}{article}{
  author={Hunt, Richard},
  author={Muckenhoupt, Benjamin},
  author={Wheeden, Richard},
  title={Weighted norm inequalities for the conjugate function and Hilbert transform},
  journal={Trans. Amer. Math. Soc.},
  volume={176},
  date={1973},
  pages={227--251},
  issn={0002-9947},
}

%\bib{MR2464252}{article}{
%   author={Hyt{\"o}nen, Tuomas},
%   title={On Petermichl's dyadic shift and the Hilbert transform},
%   language={English, with English and French summaries},
%   journal={C. R. Math. Acad. Sci. Paris},
%   volume={346},
%   date={2008},
%   number={21-22},
%   pages={1133--1136},
%   issn={1631-073X},
%   review={\MR{2464252 (2010e:42012)}},
%  %doi={10.1016/j.crma.2008.09.021},
%}


\bib{1007.4330}{article}{
  author={Hyt\"onen, Tuomas},
  title={The sharp weighted bound for general Calderon-Zygmund operators},
  eprint={http://arxiv.org/abs/1007.4330},
  date={2010},
}


\bib{1103.5229}{article}{
  author={Hyt\"onen, Tuomas},
  author={Lacey, Michael T.},
  author={Martikainen, Henri}, 
  author={Orponen, Tuomas}, 
  author={Reguera, Maria Carmen},
  author={Sawyer, Eric T.},
  author={Uriarte-Tuero, Ignacio},
  title={Weak And Strong Type Estimates for Maximal Truncations of Calder\'on-Zygmund Operators on $ A_p$ Weighted Spaces},
  eprint={http://www.arxiv.org/abs/1103.5229},
  date={2011},
}

\bib{1103.5562}{article}{
  author={Hyt\"onen, T.},
  author={P\'erez, C.},
  title={Sharp weighted bounds involving $A_\infty $},
  eprint={http://www.arxiv.org/abs/1103.5562},
  date={2011},
}



\bib{1010.0755}{article}{
  author={Hyt\"onen, T.},
  author={P{\'e}rez, Carlos},
  author={Treil, S.},
  author={Volberg, A.},
  title={Sharp weighted estimates of the dyadic shifts and $A_2$ conjecture},
  journal={ArXiv e-prints},
  date={2010}, 
  eprint={http://arxiv.org/abs/1010.0755},
}

\bib{1104.2199}{article}{
  author={Lacey, Michael T.},
  title={An $A_p$--$A_\infty$ inequality for the Hilbert Transform},
  journal={ArXiv e-prints},
  date={2011}, 
  eprint={http://arxiv.org/abs/1104.2199},
}


\bib{MR2657437}{article}{
  author={Lacey, Michael T.},
  author={Petermichl, Stefanie},
  author={Reguera, Maria Carmen},
  title={Sharp $A_2$ inequality for Haar shift operators},
  journal={Math. Ann.},
  volume={348},
  date={2010},
  number={1},
  pages={127--141},
  issn={0025-5831},
}

%\bib{0911.3437}{article}{
%  author={Lacey, Michael T.},
%  author={Sawyer, Eric T.},
%  author={Uriarte-Tuero, Ignacio},
%  title={Two Weight Inequalities for Discrete Positive Operators},
%  date={2009},
%  journal={Submitted},
%  eprint={http://www.arxiv.org/abs/0911.3437},
%}

%\bib{0807.0246}{article}{
%  author={Lacey, Michael T.},
%  author={Sawyer, Eric T.},
%  author={Uriarte-Tuero, Ignacio},
%  title={A characterization of two weight norm inequalities for maximal singular integrals with one doubling measure},
%  date={2008},
%  journal={ A\&PDE, to appear},
%  eprint={http://arxiv.org/abs/0805.0246},
%}
%
%\bib{0911.3920}{article}{
%  author={Lacey, Michael T.},
%  author={Sawyer, Eric T.},
%  author={Uriarte-Tuero, Ignacio},
%  title={Two Weight Inequalities for Maximal Truncations of Dyadic Calder\'on-Zygmund Operators},
%  date={2009},
%  journal={Submitted},
%  eprint={http://www.arxiv.org/abs/0911.3920},
%}


\bib{MR2721744}{article}{
   author={Lerner, Andrei K.},
   title={A pointwise estimate for the local sharp maximal function with
   applications to singular integrals},
   journal={Bull. Lond. Math. Soc.},
   volume={42},
   date={2010},
   number={5},
   pages={843--856},
   issn={0024-6093},
   review={\MR{2721744}},
  %doi={10.1112/blms/bdq042},
}


%\bib{1005.1422}{article}{
%  author={Lerner, Andrei K.},
%  title={Sharp weighted norm inequalities for Littlewood-Paley operators and singular integrals},
%  date={2010},
%  eprint={http://arxiv.org/abs/1005.1422},
%}

%\bib{lerner}{article}{
%  author={Lerner, Andrei K.},
%  title={On some weighted norm inequalities for Littlewood-Paley operators},
%  journal={Illinois J. Math.},
%  volume={52},
%  date={2007},
%  number={2},
%  pages={653--666},
%}





%
%\bib{MR1756958}{article}{
%   author={Petermichl, Stefanie},
%   title={Dyadic shifts and a logarithmic estimate for Hankel operators with
%   matrix symbol},
%   language={English, with English and French summaries},
%   journal={C. R. Acad. Sci. Paris S\'er. I Math.},
%   volume={330},
%   date={2000},
%   number={6},
%   pages={455--460},
%   issn={0764-4442},
%   review={\MR{1756958 (2000m:42016)}},
%}

\bib{MR2354322}{article}{
  author={Petermichl, Stefanie},
  title={The sharp bound for the Hilbert transform on weighted Lebesgue spaces in terms of the classical $A\sb p$ characteristic},
  journal={Amer. J. Math.},
  volume={129},
  date={2007},
  number={5},
  pages={1355--1375},
  issn={0002-9327},
}


%\bib{MR1964822}{article}{
%   author={Petermichl, S.},
%   author={Treil, S.},
%   author={Volberg, A.},
%   title={Why the Riesz transforms are averages of the dyadic shifts?},
%   booktitle={Proceedings of the 6th International Conference on Harmonic
%   Analysis and Partial Differential Equations (El Escorial, 2000)},
%   journal={Publ. Mat.},
%   date={2002},
%   number={Vol. Extra},
%   pages={209--228},
%   issn={0214-1493},
%   review={\MR{1964822 (2003m:42028)}},
%}

\bib{1103.5347}{article}{
  author={Reznikov, Alexander},
  author={Treil, Sergei},
  author={Volberg, Alexander},
  title={A sharp weighted estimate of dyadic shifts of complexity 0 and 1}, 
  journal={ArXiv e-prints},
  date={2011}, 
  eprint={http://arxiv.org/abs/1103.5347},
}




%\bib{MR676801}{article}{
%   author={Sawyer, Eric T.},
%   title={A characterization of a two-weight norm inequality for maximal
%   operators},
%   journal={Studia Math.},
%   volume={75},
%   date={1982},
%   number={1},
%   pages={1--11},
%   issn={0039-3223},
%   review={\MR{676801 (84i:42032)}},
%}
%
%\bib{MR719674}{article}{
%  author={Sawyer, Eric},
%  title={A two weight weak type inequality for fractional integrals},
%  journal={Trans. Amer. Math. Soc.},
%  volume={281},
%  date={1984},
%  number={1},
%  pages={339--345},
%  issn={0002-9947},
%}
%
%\bib{MR930072}{article}{
%  author={Sawyer, Eric T.},
%  title={A characterization of two weight norm inequalities for fractional and Poisson integrals},
%  journal={Trans. Amer. Math. Soc.},
%  volume={308},
%  date={1988},
%  number={2},
%  pages={533--545},
%  issn={0002-9947},
%}

\bib{1105.2252}{article}{
   author={{Treil}, S.},
   title={Sharp $A_2$ estimates of Haar shifts via Bellman function},
   year={2011},
   eprint={http://arxiv.org/abs/1105.2252},
}



%\bib{MR883661}{article}{
%  author={Wilson, J. Michael},
%  title={Weighted inequalities for the dyadic square function without dyadic {$A_\infty $}},
%  journal={Duke Math. J.},
%  fjournal={Duke Mathematical Journal},
%  volume={55},
%  year={1987},
%  number={1},
%  pages={19--50},
%  issn={0012-7094},
%  coden={DUMJAO},
%  mrclass={42B25},
%  mrnumber={883661 (88d:42034)},
%  mrreviewer={B. Muckenhoupt},
%  %doi={10.1215/S0012-7094-87-05502-5},
%  url={http://dx.doi.org/10.1215/S0012-7094-87-05502-5},
%}
%
%
%\bib{MR2680056}{article}{
%   author={Vagharshakyan, Armen},
%   title={Recovering singular integrals from Haar shifts},
%   journal={Proc. Amer. Math. Soc.},
%   volume={138},
%   date={2010},
%   number={12},
%   pages={4303--4309},
%   issn={0002-9939},
%   review={\MR{2680056}},
%  %doi={10.1090/S0002-9939-2010-10426-4},
%}


\end{biblist}
\end{bibsection}
\end{document}